\documentclass[12pt,reqno]{article}
\usepackage{amsmath,amsthm,amsfonts,amssymb,amscd}
\usepackage[dvips]{graphicx}
\usepackage{psfrag}

\input{psfig}
\theoremstyle{plain}
\newtheorem{thm}{Theorem}[section]
\newtheorem{rk}[thm]{Remark}

\newtheorem{clly}[thm]{Corollary}
\newtheorem{lemma}[thm]{Lemma}
\newtheorem{defi}[thm]{Definition}
\newtheorem{nota}[thm]{Notation}
\newtheorem{maintheorem}{Theorem}

\title{The explosion of singular hyperbolic attractors}
\author{C. A. Morales
        \thanks{
This work is partially supported by CNPq, FAPERJ and PRONEX-Dyn. Systems/Brazil.}}

\begin{document}
\maketitle
\begin{abstract}
A {\em singular hyperbolic attractor} for
flows is a partially hyperbolic attractor
with singularities (hyperbolic ones)
and volume expanding central direction \cite{mpp1}.
The geometric Lorenz attractor \cite{gw}
is an example of a
singular hyperbolic attractor.
In this paper we study the perturbations
of singular hyperbolic attractors
for three-dimensional flows.
It is proved that any attractor
obtained from such perturbations contains a singularity.
So,
there is an upper bound for the number
of attractors obtained from such perturbations.
Furthermore, every three-dimensional flow
$C^r$ close to one exhibiting a
singular hyperbolic attractor has a singularity
non isolated in the non wandering set.
We also give sufficient
conditions for a singularity of a three-dimensional
flow to be stably non isolated in the nonwandering set.
These results generalize well known properties of the Lorenz attractor.
\end{abstract}

\section{Introduction}

In \cite{mpp1}
it was proved that
$C^1$ robust attractors with
singularities for three-dimensional $C^1$ flows
are singular hyperbolic.
A natural question is whether
singular hyperbolic attractors
for three-dimensional $C^1$ flows are $C^1$ robust, but
the answer is negative
by \cite{mp}.
This fact motivates the
study the perturbations of singular hyperbolic attractors
for three-dimensional flows.
We shall prove that
any attractor obtained from such perturbations
is singular (i.e. contains a singularity).
Consequently there is an
upper bound for the number of attractors
obtained from such perturbations.
Furthermore, every three-dimensional flow
$C^r$ close to one exhibiting
singular hyperbolic attractors has a singularity
non isolated in the non wandering set.
We also give a sufficient condition for a singularity
of a three-dimensional flow to be stably
non isolated in the nonwandering set.
These results generalize
well known properties of the
geometric Lorenz attractor \cite{gw}.

Let us describe our results in a precise way.
Hereafter $M$ is a closed
$3$-manifolds and
${\cal X}^r(M)$ denotes the space
of $C^r$ flows in $M$ endowed with the
$C^r$ topology, $r\geq 1$.
Every flow $X=X_t\in {\cal X}^r(M)$ will be assumed to be orientable and we still denote by $X$ its correponding vector field.
Hereafter we fix a flow $X\in {\cal X}^r(M)$.
A point $p$ is
{\em nonwandering} if
$\forall U$ neighborhood
of $p$ and $\forall T>0$
$\exists t>T$ such that $X_t(U)\cap U\neq\emptyset$.
The {\em nonwandering set} of $X$ is the set of
nonwandering points
of $X$.
The $\omega$-limit set of $p$ is the set
$\omega_X(p)$ of $q$ such that $q=\lim_{n\infty}X_{t_n}(p)$ for some
sequence $t_n\to\infty$.
An invariant set is {\em transitive} if
it is the $\omega$-limit set of one of its points.
A compact invariant set is :
{\em singular} if it has a singularity and
{\em attracting}
if it is $\cap_{t>0}X_t(U)$
for some compact neighborhood $U$ of it satisfying $X_t(U)\subset U$,
$\forall t>0$ ($U$ is called
isolating block).
An {\em attractor} is a transitive attracting set.
An attractor is : {\em nontrivial} if it is
not a single closed orbit and
{\em $C^r$-robust} if there is an isolating block $U$ of it such that $\cap_{t>0}Y_t(U)$ is both
transitive and nontrivial
for every flow $Y$ $C^r$ close to $X$.
A compact invariant set is {\em hyperbolic}
if it exhibits a tangent bundle decomposition
$E^s\oplus E^X\oplus E^u$ such that $E^s$ is contracting,
$E^u$ is expanding and $E^X$ is the direction of $X$.
A closed orbit is hyperbolic if its orbit
is hyperbolic as a compact invariant set.
Nontrivial hyperbolic attractors
are often called {\em hyperbolic strange attractors}
(this is equivalent
to $E^s\neq 0$ for the corresponding hyperbolic splitting).
A compact, singular, invariant set
of $X$ is {\em singular hyperbolic}
if all its singularities
are hyperbolic and, in addition,
it exhibits a continuous invariant splitting
$E^s\oplus E^c$ such that
$E^s$ is contracting, $E^s$ dominates
$E^c$ and
$E^c$ is volume expanding
(i.e. the jacobian of the derivative of $DX_t$ along $E^c$ growths exponentially as $t\to \infty$).
The singular hyperbolic sets
belongs to the category
of partially hyperbolic sets.
A {\em singular hyperbolic attractor}
is a singular hyperbolic set which is also
an attractor.
The most representative example of a singular
hyperbolic attractor is the {\em geometric
Lorenz attractor} \cite{gw}.
Our first result is the following.

\begin{maintheorem}
\label{thB}
For every
singular hyperbolic attractor of $X$
there is a neighborhood $U$ of it such that
every attractor in $U$ of every flow $C^r$ close to $X$
is singular.
\end{maintheorem}

Let us present some applications of Theorem \ref{thB}
for the study of the perturbations of singular hyperbolic attractors in dimension three.
First we prove the existence of an upper bound
for the number of attractors obtained by
perturbing a singular hyperbolic one. 

\begin{clly}
\label{bound}
For every singular hyperbolic attractor of $X$
there are a neighborhood $U$ of it and $n
\in I\!\! N^*$ such that every flow $C^r$ close to $X$
has at most $n$ attractors in $U$.
\end{clly}

Second we study
the existence of
flows exhibiting in a stably way singularities
non isolated in the nonwadering set.
By definition
a singularity $\sigma$ of a flow $X$ is {\em isolated
in the nonwandering set}
if the set
$\Omega(X)\setminus \{\sigma\}$ is
not closed in $M$.
The existence of singularities
non isolated in the nonwandering set
is a well known obstruction for hyperbolicity
(this fact was explored in \cite[Theorem 3]{n}).
Observe that
every flow close to one exhibiting
a geometric Lorenz attractor has a singularity
non isolated in the nonwandering set.
In general we observe that every flow $C^r$ close to one exhibiting a $C^r$ robust
singular attractor has a singularity
non isolated in the nonwandering set.
As
$C^1$ robust singular attractors
on closed $3$-manifolds are singular hyperbolic
(but not conversely)
it is natural to ask
whether  the conclusion of the last observation
holds for singular hyperbolic attractors instead of
robust singular ones.
The answer is positive by the following result.

\begin{clly}
\label{thA}
Every flow in ${\cal X}^r(M)$ which is $C^r$ close
to one exhibiting singular hyperbolic attractors has a singularity non isolated in the nonwandering set.
\end{clly}

Third we
proves the conclusion above for {\em generic}
three-dimensional $C^1$ flows exhibiting nontrivial singular attractors.
Recall that a subset of ${\cal X}^r(M)$ is
{\em residual} whenever it
contains a countable intersection of
open-dense subsets of ${\cal X}^r(M)$.

\begin{clly}
\label{c1}
Every flow in a residual subset of
${\cal X}^1(M)$ exhibiting a nontrivial singular attractor has a singularity non isolated in the
nonwandering set.
\end{clly}

Finally we discuss the existence of singularities
stably non isolated in the nonwandering set.
By definition a hyperbolic singularity
$\sigma$ of $X$ is {\em $C^r$ stably
non isolated in the nonwandering set}
if
$\Omega(Y)\setminus \{\sigma(Y)$
is not closed,
$\forall Y$ $C^r$ close to $X$, where
$\sigma(Y)$ denotes the continuation of $\sigma$ for
$Y$ close to $X$ \cite{dmp}.
An example of a singularity
stably non isolated in the nonwanering
set is the geometric Lorenz attractor's one.
The following result which easily follows from
Theorem \ref{thB} and Corollary \ref{thA}
gives a sufficient condition for a singularity
to be stably non isolated in the nonwandering set.
See \cite{mp3} for a sort of converse when $r=1$.

\begin{clly}
\label{Coro}
If $\sigma$ is the {\em unique} singularity
of a singular hyperbolic attractor
of a flow in ${\cal X}^r(M)$,
then
$\sigma$ is $C^r$ stably non isolated in the nonwandering set.
\end{clly}

The proof of Theorem \ref{thB}
is as follows.
By contradiction we suppose that there are
a three-dimensional $C^r$ flow
$X$, with a singular hyperbolic attractor
$\Lambda$, and a sequence $X^n$ of $C^r$ flows
converging to $X$ such that
every $X^n$ exhibits a non singular attractor $A^n$
in a fixed isolating block
$U$ of $\Lambda$. If $n$ is large then
$U$ is attracting for $X^n$ too.
It will follow that $A^n$
is a hyperbolic strange attractor
of $X^n$.
In particular, as $dim(M)=3$,
the strong unstable manifolds of $A^n$ are one-dimensional.
We have two cases,
namely the sequence $A^n$ either accumulates
some singularity $\sigma$ of $X$ or does not.
In the later case
the invariant manifold
theory  and the stability of hyperbolic sets
will yield a contradiction.
In the former case we
shall prove that
$A^n$ intersects the stable manifold
of $\sigma(X^n)$, the continuation
of $\sigma$. For this we analyze the
relative position of the stable
and central direction of the umperturbed
attractor close to $\sigma$.
Such directions depends continuosly
on the flow.
This is used to prove
$\sigma(X^n)\in A^n$ a contradiction
since $A^n$ is strange.
We derive the corollaries
\ref{bound}, \ref{thA}, \ref{c1} from Theorem \ref{thB}
in the last section.

\section{Preliminars}

In this section we establish some definitions and preliminar results.
Throughout $M$ denotes a closed $3$-manifold
and $X$ denotes a $C^r$ flow in $M$,
$r\geq 1$.
We denote by
$Sing(X)$ the set of singularities
of $Y$. Given $A\subset M$
we denote by $Cl(A)$ the closure of $A$.
If $\delta>0$ we
denote
$$
B_\delta(A)=\{x\in M:d(x,A)<\delta\},
$$
where $d(\cdot,\cdot)$ is the metric in
$M$.

Given $p\in M$ and $\epsilon>0$
we define
$$
W^{ss}_X(p)=
\{x:d(X_t(x),X_t(p))\to 0, t\to \infty\},
$$
$$
W^{uu}_X(p)=
\{x:d(X_t(x),X_t(p))\to 0, t\to -\infty\},
$$
$$
W^{ss}_X(p,\epsilon)=
\{x:d(X_t(x),X_t(p))\leq\epsilon, \forall t\geq 0\}
$$
and
$$
W^{uu}_X(p,\epsilon)=
\{x:d(X_t(x),X_t(p))\leq \epsilon, \forall
t\leq 0\}.
$$
These sets are called respectively
the stable, unstable, local stable and local unstable
set of $p$.
We also define
$$
W^{s}_X(p)=\cup_{t\in I\!\! R}W^{ss}_X(X_t(p))
$$
and
$$
W^{u}_X(p)=\cup_{t\in I\!\! R}W^{uu}_X(X_t(p)).
$$

A compact invariant set $H$ of $X$ is {\em
hyperbolic}
if there are a continuous tangent bundle
splitting $E^s\oplus E^X\oplus E^u$ over
$H$ and positive constants
$C,\lambda$ such that

\begin{itemize}
\item $E^X$ is the flow's
direction.
\item $E^s$ is contracting, i.e.
$$
\mid\mid DX_t(x)/E^s_x\mid\mid
\leq Ce^{-\lambda t}$$
for all $x\in H$ and $t>0$.
\item $E^u$ is expanding, i.e.
$$
\mid\mid DX_{-t}(x)/E^u_x\mid\mid
\leq Ce^{-\lambda t},
$$
for all $x\in H$ and $t> 0$.
\end{itemize}

If $H$ is a hyperbolic set of $X$
and $p\in H$ then
$W^{ss}_X(p)$ and $W^{s}_X(p)$
are $C^r$-submanifolds
of $M$. This is the Stable Manifold Theorem
\cite{hps}.
Similarly
$W^{uu}_X(p)$, $W^u_X(p)$,
$W^{ss}_X(p,\epsilon)$ and $W^{uu}_X(p,\epsilon)$
are $C^r$ submanifolds.
In this case
$W^{ss}_X(p)$ and $W^{uu}_X(p)$ are called
respectively
the {\em strong stable} and the {\em
strong unstable} manifold
of $p$.
Note that $dim(W^{ss}_X(p))=dim(E^s)$
and $dim(W^{uu}_X(p)=dim(E^u)$.
If $O$ is a closed orbit of $X$ then
we denote by $O(Y)$ the continuation of $O$
for $Y$ $C^r$ close to $X$ \cite{dmp}.

\begin{defi}
\label{d2}
A compact invariant set
$\Lambda$ of $X$
is {\em partially hyperbolic}
if there are a continuous tangent bundle
splitting $E^s\oplus E^c$ over
$\Lambda$ and positive constants $C,\lambda$
such that
\begin{enumerate}
\item
$E^s$ is expanding, i.e.
$$
\mid\mid DX_t(x)/E_x^s\mid\mid
\leq Ce^{-\lambda t},
$$
for all $x\in \Lambda$ and $t>0$.
\item
$E^s$ dominates $E^c$, i.e.
$$
\mid\mid DX_t(x)/E_x^s\mid\mid\cdot
\mid\mid DX_{-t}(x)/E^c_x\mid\mid
\leq Ce^{-\lambda t},
$$
for all $x\in \Lambda$ and $t>0$.
\end{enumerate}
The subbundle $E^c$ above is
the {\em central direction} of $\Lambda$.
We say that $E^c$ is {\em volume expanding} if for every $x\in \Lambda$ we have
$$
\mid det DX_t(x)/E^c_x\mid
\geq C^{-1}e^{\lambda t},
$$
where $det DX_t(x)/E^c_x$ denotes
the jacobian of $DX_t(x)$ along $E^c_x$.
\end{defi}

\begin{defi}
\label{shs}(\cite{mpp2})
A {\em singular hyperbolic set}
is a partially hyperbolic set having
singularities (all of them hyperbolic)
and volume expanding central direction.
A {\em singular hyperbolic attractor}
is an attractor which is also
a singular hyperbolic set.
\end{defi}

\begin{rk}
\label{r3}
\hspace{1pt}
\begin{enumerate}
\item
A singular hyperbolic attractor
cannot be a hyperbolic set.
The most well known example of singular hyperbolic
attractor is the geometric Lorenz attractor \cite{gw}.
\item
Let $T$ be a transitive, nontrivial, singular hyperbolic set
of $X\in {\cal X}^r(M)$. Denote by
$E^s\oplus E^c$ the singular hyperbolic splitting
of $T$. Then
$dim(E^s)=1$, $dim(E^c)=2$ and the subspace
$E^X$ generated by
$X$ in the tangent space
is contained in
$E^c$.
\end{enumerate}
\end{rk}

An useful property of singular hyperbolic sets
is given below.

\begin{lemma}(\cite{mpp2})
\label{l1}
Let $\Lambda$ be a singular hyperbolic set of
a $C^r$ flow
$X$ in $M$, $r\geq 1$.
There is a neighborhood $U$ of $\Lambda$
such that if $Y$ is $C^r$ close to $X$,
then every nonempty, compact, non singular, invariant set
of $Y$ in $U$ is
hyperbolic.
\end{lemma}

\begin{lemma}
\label{l1'}
Let $\Lambda$ be a singular hyperbolic attractor
of a $C^r$ flow $X$ in $M$, $r\geq 1$.
If $O\subset \Lambda$ is a periodic orbit
of $X$,
then $O$ is a hyperbolic saddle-type periodic orbit.
Henceforth $W^{uu}_X(p)$ is a one-dimensional
submanifold.
In addition, for every $p\in O$
the set
$$
\{q\in W^{uu}_X(p):
\Lambda=\omega_X(q)\}
$$
is dense in $W^{uu}_X(p)$.
\end{lemma}

\begin{proof}
First observe that the contracting direction
$E^s$ of $\Lambda$ is integrable \cite{hps}.
Denote by $W^{ss}_X(p)$ the leave of the
resulting foliation passing
throught $p\in \Lambda$.
Note that the forward orbit
of every $x\in W^{ss}_X(p)$
is asymptotic to the forward orbit
of $p$.
If follows from
Lemma \ref{l1} that $O$ is hyperbolic.
That $O$ is saddle-type follows
because $E^s$ is contracting and $E^c$ volume expanding
(see definitions \ref{d2} and \ref{shs}).
In particular $dim(W^{uu}_X(p))=1$ for all
$p\in O$.
Let $I$
be an open interval in $W^{uu}_X(p)$.
It follows that
$$
B=\cup_{0\leq t\leq 1} X_t(I)
$$
is a two-dimensional manifold of $M$.
So,
$$
B'=\cup_{x\in B} W^{ss}_X(x)
$$
contains an open set $V$ with
$B\cap V\neq \emptyset$.
Clearly $B\cap V\subset \Lambda$ since
$\Lambda$ is an attractor,
$B\subset W^u_X(p)$ and $p\in O$.
Let $q\in \Lambda$ such that
$\Lambda=\omega_X(q)$.
Then the forward orbit
of $q$ intersects $V$ and so it intersects
$B'$ too.
It follows from the definition of
$B'$ that the positive orbit
of $q$ is asymptotic to
the forward orbit of some
$q'\in B$.
In particular,
$\Lambda=\omega_X(q)=\omega_X(q')$.
This proves that
$
\{q\in W^{uu}_X(p):
\Lambda=\omega_X(q)\}
$ is dense in $\Lambda$ as desired.
\end{proof}

We shall use the following definition.
A singularity
$\gamma$ of $X$ is {\em Lorenz-like}
if its eigenvalues 
$\lambda_1,\lambda_2,\lambda_3$ are real and
satisfy
$$
\lambda_2<\lambda_3<0<-\lambda_3<\lambda_1.
$$

By the Invariant Manifold Theory \cite{hps}
it follows that the eigenspace
associated to $\lambda_2$
is tangent to a one-dimensional invariant 
manifold of $X$ denoted by $W^{ss}_X(\gamma)$.
As usual we denote by $W^s_{loc}(\gamma)$ a small ball in $W^s_X(\gamma)$
centered at $\gamma$. 
Similarly we define $W^u_{loc}(\gamma)$.
All these invariant manifolds persist under
small $C^r$ perturbations $Y$ of $X$.

We consider Lorenz-like singularities due to the following result.

\begin{thm}(\cite{mpp2})
\label{t1}
If $\Lambda$ is a singular hyperbolic
attractor of $X$ then
every $\gamma\in Sing(X)\cap \Lambda$ is
is Lorenz-like
and satisfies
$$
\Lambda\cap W^{ss}_X(\gamma)=\{\gamma\}.
$$
\end{thm}

By \cite[Hartman Grobman Theorem]{dmp}
we can describe
the flow of $X$ around a Lorenz-like singularity $\gamma$ (see Figure 1).

\begin{figure}[htv] 
\centerline{
\psfig{figure=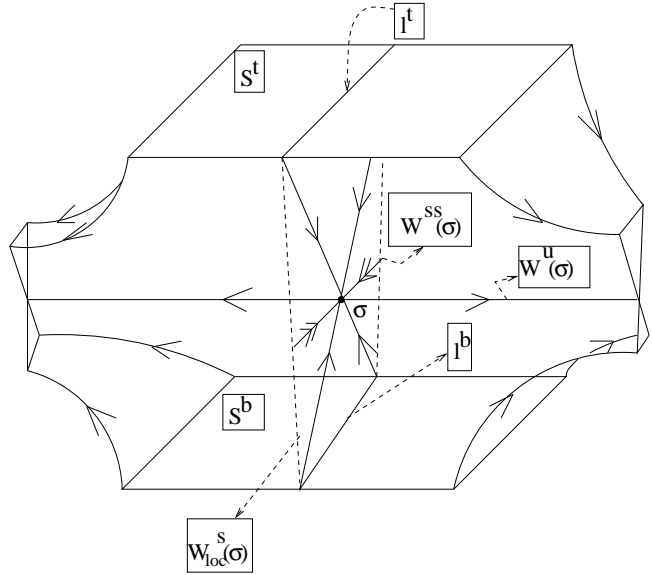,height=3in}}
\caption{\label{f.1} The flow around $\sigma$.}
\end{figure}

Indeed, $\gamma$ has a two-dimensional stable manifold $W^s_X(\gamma)$, 
a strong stable manifold $W^{ss}_X(\gamma)$, and an 
unstable manifold $W^u_X(\gamma)$ both one-dimensional.
It follows that $W^{ss}_X(\gamma)$ separates $W^s_{loc}(\gamma)$
in two connected components : the top and the bottom ones.

\begin{nota}
\label{n1}
We denote by $S^t$ and $l^t$ the
cross-section and the curve as in Figure 
Similarly we denote $S^b$
and $l^b$.
Topologically $S^t$ is a rectangle
and so the the boundary of $S^t$ is formed by four curves,
two of them being transverse to $l^t$.
We denote by $\partial^h S^t$ the union of such curves. 
Similarly we define $\partial^h S^b$.
The interior of $S^t$ as a submanifold of $M$
is denoted by $int(S^t)$.
Similarly we define $int(S^b)$.
\end{nota}

\begin{rk}
\label{r4}
As mentioned
in Notation \ref{n1}, $S^t$ and $S^b$ are
rectangles, i.e.
diffeomorphic to $[-1,1]\times [-1,1]$
(say). Note that $l^t$ and $l^b$ are contained 
in $W^s_X(\gamma)\setminus W^{ss}_X(\gamma)$.
Clearly $S^t,S^b$ can be chosen to be
close to $\gamma$.
The positive flow lines of $X$ starting at $S^t\cup S^b\setminus(l^t\cup l^b)$
exit a small neighborhood of $\gamma$ passing through the cusp region
as indicated in Figure
The positive orbits starting at $l^t\cup l^b$ goes directly to $\gamma$.
\end{rk}


\section{Proof of Theorem \ref{thB}}

Let $X$ and $\Lambda$ be as in the hypothesis
of the theorem.
Fix the neighborhood
$U$ of $\Lambda$ obtained
in the conclusion of Lemma \ref{l1}.
If $U'$ is an isolating block
of $\Lambda$ then there is $T>0$ such that $X_T(U')
\subset U$.
Clearly $X_T(U')$ is also
an isolating block of $\Lambda$.
In particular, we can assume
that $U'\subset U$.
We can further assume that $U=U'$ without
loss of generality.

By contradiction suppose that
the conclusion of the theorem is not true for
$X,\Lambda,U$. Then
there is a sequence of flows
$X^n$ converging to
$X$ in the $C^r$ topology such that
every $X^n$ has an
attractor
$A^n\subset U$ which is not singular hyperbolic.
As $X^n\to X$
we can assume that $X^n$ is $C^r$ close to
$X$
for all $n$.
Then, by lemmas \ref{l1} and \ref{l1'},
$A^n$ is a hyperbolic strange attractor
of $X^n$.

We shall obtain a contradiction by using the
attractors $A^n$.
For this we separate the proof in two cases :

\begin{description}
\item{Case I:}
$Sing(X)\cap Cl\left(\cup_{n\in I\!\! N}A^n\right)
=\emptyset$.
\item{Case II:}
$Sing(X)\cap Cl\left(\cup_{n\in I\!\! N}A^n\right)
\neq\emptyset$.
\end{description}

\noindent
{\flushleft{\bf Proof in Case I: }}
In this case there is
$\delta>0$ such that
\begin{equation}
\label{*}
B_\delta(Sing(X))\cap\left(\cup_{n\in I\!\! N}A^n\right)=\emptyset.
\end{equation}

Fix a compact neighborhood
$V\subset U$ of $\Lambda$ such that
$$
\Lambda=\cap_{t\geq 0}X_t(V).
$$
As $U$ is an isolating block of
$\Lambda$ can choose $V$ so close to $U$
to ensure
$$
A^n\subset V
$$ for all $n$.
Define
$$
H=
\cap_{t\geq 0}X_t\left(V\setminus B_{\delta/2}(Sing(X))
\right).
$$
Obviously $Sing(X)\cap H=\emptyset$.

Clearly $H\neq\emptyset$.
Indeed, for all $n$ we choose $x^n\in A^n$.
By passing to a subsequence if necessary
we can assume that
$x^n$ converges to some
$x\in M$.
Clearly
$x\in H$
for, otherwise,
the $X^n$-orbit of $x^n$ would intersect
$B_\delta(Sing(X))$
by \cite[Tubular Flow-Box Theorem]{dmp}.
This is impossible by (\ref{*}) proving
$H\neq\emptyset$ as desired.

Note that $H$ is
compact since $V$ is.
It follows that $H$ is a nonempty compact
invariant set of $X$. As
$Sing(X)\cap H=\emptyset$ it follows that $H$
is hyperbolic by Lemma \ref{l1}.
Denote by $E^s\oplus E^X\oplus E^u$ the corresponding hyperbolic splitting.

By the stability of hyperbolic sets we can fix a neighborhood $W$ of
$H$ and $\epsilon>0$ such that
if $Y$ is a flow
$C^r$ close to $X$
and $H_Y$ is a compact invariant set
of $Y$ in $W$ then :

\begin{description}
\item{(H1).}
$H_Y$ is hyperbolic and its
hyperbolic splitting
$E^{s,Y}\oplus E^Y\oplus
E^{u,Y}$
satisfies
$$
dim(E^u)=dim(E^{u,Y}),\,\,\,\,
dim(E^s)=dim(E^{s,Y}).
$$
\item{(H2).}
The manifolds
$W^{uu}_Y(x,\epsilon)$, $x\in H_Y$, are one-dimensional
and have uniform size $\epsilon$.
\end{description}

As $X^n\to X$ we have that
$$
\cap_{t\in I\!\! R}X^n_t(V\setminus B_{\delta/2}(Sing(X))
\subset W
$$
for all $n$ large.
But
$$
A^n\subset
V\setminus B_{\delta/2}(Sing(X))
$$
for all $n$
by (\ref{*}).
As $A^n$ is $X^n$-invariant we conclude
that
$$
A^n\subset W
$$
for all $n$ large.

By (H2)
we have that
$W^{uu}_{X^n}(x^n,\epsilon)$
is one-dimensional and
has uniform size $\epsilon$,
$\forall x^n\in A^n$ and $\forall n$.

Choose $x^n\in A^n$
so that $x^n$ converges to some $x\in M$.
As observed above we have that
$x\in H$.
Note that
the tangent vectors
of the curve $W^{uu}_{X^n}(x^n,\epsilon)$
at $c\in W^{uu}_{X^n}(x^n,\epsilon)$
is in $E^{u,X^n}_c$
($E^{u,X^n}$ is defined by (H1) as $X^n\to X$).

As $X^n\to X$ we have that
the angle between the directions
$E^{u,X^n}$ and $E^u$ goes to zero
as $n\to\infty$.
Henceforth
the manifolds
$W^{uu}_X(x,\epsilon)$ and
$W^{uu}_{X^n}(x^n,\epsilon)$
are almost paralell as $n\to\infty$.
As $x^n\to x$ we conclude that 
$$
W^{uu}_{X^n}(x^n,\epsilon)
\to W^{uu}_X(x,\epsilon)
$$
in the sence of $C^1$ submanifolds
\cite{pt}.

Fix an open interval
$I\subset
W^{uu}_X(x,\epsilon)$ containing $x$.
By Lemma \ref{l1'}, as $\Lambda\cap Sing(X)\neq
\emptyset$, we have that
there are $q\in I$ and $T>0$ such that
$$
X_T(q)\in B_{\delta/5}(Sing(X)).
$$
By \cite[Tubular Flow Box Theorem]{dmp} there is $V_q$ open containing $q$ such that
$$
X_T(V_q)\subset B_{\delta/5}(Sing(X)).
$$
As $X^n\to X$
we have
\begin{equation}
\label{intersection}
X^n_T(V_q)\subset B_{\delta/4}(Sing(X))
\end{equation}
for all $n$ large.

But $
W^{uu}_{X^n}(x^n,\epsilon)
\to W^{uu}_X(x,\epsilon)
$, $q\in I\subset W^{uu}_X(p,\epsilon)$,
$q\in V_q$ and $V_q$ is open.
So,
$$
W^{uu}_{X^n}(x^n,\epsilon)
\cap V_q\neq\emptyset
$$
for all $n$ large.
Applying (\ref{intersection})
to $X^n$ for $n$ large we have
$$
X^n_T(W^{uu}_{X^n}(x^n,\epsilon))\cap
B_{\delta/4}(Sing(X))\neq\emptyset.
$$
As $W^{uu}_{X^n}(x^n,\epsilon)\subset W^u_{X^n}(x^n)$
the invariance of $ W^u_{X^n}(x^n)$
implies
$$
W^u_{X^n}(x^n)\cap B_{\delta/2}(Sing(X))
\neq\emptyset.
$$
Observe that $W^u_X(x^n)\subset A^n$
since
$x^n\in A^n$ and $A^n$ is an attractor.
We conclude that
$$
A^n\cap B_{\delta}(Sing(X))\neq\emptyset.
$$
This contradicts
(\ref{*}) and
the proof follows in Case I.

\qed

\noindent
{\flushleft{\bf Proof in Case II: }}
In this case there is $\sigma\in Sing(X)\cap \Lambda$
so that
$$
\sigma\in
Cl\left(\cup_{n\in I\!\! N}A^n\right).
$$ 

By Theorem \ref{t1} we have that
$\gamma=\sigma$ is Lorenz-like
and satisfies
$$
\Lambda\cap W^{ss}_X(\sigma)=\{\sigma\}.
$$

As $\sigma$ is Lorenz-like,
the cross-sections
$S^t,S^b$ at Notation \ref{n1} are well defined.
The last equality implies
that $S^t,S^b$
can be chosen so that
$$
\Lambda\cap\left(\partial^hS^t\cup\partial^hS^b\right)
=\emptyset.
$$

As $X^n\to X$ we have that
$S^t,S^b$ are
cross-sections of $X^n$ too.
By using the Implicit Function
Theorem we can assume that $\sigma(X^n)=\sigma$
and
\begin{equation}
\label{l}
l^t\cup l^b\subset W^s_{X^n}(\sigma)
\end{equation}
for all $n$.
We shall denote by $TB=T(B)$ the tangent space at $B$.

Now,
the one-dimensional subbundle $E^s$ of $\Lambda$
extends to a contracting invariant
subbundle on the hole $U$.
Take a continuous (not necessarily invariant)
extension of $E^c$
to $U$.
As usual we still denoted the above
extensions by $E^s\oplus E^c$.

By the invariant manifold
theory [HPS] it follows
that the splitting
$E^s\oplus E^c$ persists by small
perturbations of $X$.
More precisely, for all $n$ large
the flow $X^n$ has an splitting
$E^{s,n}\oplus E^{c,n}$ over $U$
such that
$E^{s,n}$ is invariant contracting,
$E^{s,n}\to E^s$ and
$E^{c,n}\to E^c$ as $n\to\infty$.
Choosing $S^t\cup S^b$ close to $\sigma$ if
necessary (Remark \ref{r4})
we can assume
that $S^t\cup S^b\subset U$.
In particular, $E^{s,n}\oplus E^{c,n}$
is defined
in $S^t\cup S^b$ for
all $n$ large. In what follows we denote
by $E^Y$ the subbundle in $TM$ generated by
a flow $Y$ in $M$.

The dominance condition
Definition \ref{d2}-(2)
and \cite[Proposition 2.2]{d}
imply that for $i=t,b$ one has
$$
T_xS^i\cap \left(E^s_x\oplus E^X_x\right)
=T_xl^i,
$$
for all $x\in l^i$.

Denote by $\angle(E,F)$ the angle between
two linear subspaces.
The last equality implies
that there is $\rho>0$ such that
$$
\angle(T_xS^i\cap E_x^c,T_xl^i)>\rho,
$$
for all $x\in l^i$
($i=t,b$).
But $E^{c,n}\to E^c$ as $n\to\infty$.
So for all $n$ large we have
\begin{equation}
\label{angle}
\angle(T_xS^i\cap E_x^{c,n},T_xl^i)>\frac{\rho}{2},
\end{equation}
for all $x\in l^i$
($i=t,b$).

Fix a coordinate system $(x^i,y^i)$
in $S^i$ ($i=t,b$) such that
$$
S^i=[-1,1]\times[-1,1],
\,\,\,\,\,\,\,\,l^i=\{0\}\times[-1,1]
$$
with respect to $(x^i,y^i)$ (Remark \ref{r4}).

Denote by
$$
\Pi^i:S^i\to [-1,1]
$$
the projection $\Pi^i(x^i,y^i)=x^i$.
Given $\Delta>0$ we define
$$
l^i(\Delta)=[-\Delta,\Delta]\times [-1,1].
$$
The horizontal boundary of $l^i(\Delta)$ is defined by
$$
\partial^hl^i(\Delta)=
[-\Delta,\Delta]\times \{\pm 1\}.
$$

The inequality (\ref{angle}) and the continuity
of $E^{c,n}$ imply
that there is $\Delta_0>0$ such that for all $n$ large
the line field $F^n$ in $l^i(\Delta_0)$ defined
by
$$
F^n_x=
T_xS^i\cap E^{c,n}_x, \,\,\,\,x\in l^i(\Delta_0)
$$
is transverse to $\Pi^i$.

As $\sigma\in Cl\left(\cup_{n\in I\!\! N}A^n\right)$,
there is $n_0$ such that
$A^{n_0}$ intersects $int(l^t(\Delta_0))$ or
$int(l^b(\Delta_0))$.
For simplicity we
denote $Z=X^{n_0}$, $A=A^{n_0}$, $F=F^{n_0}$.
In particular $A$ is a hyperbolic strange
attractor of $Z$.

We can assume $A\cap int(l^t(\Delta_0))\neq\emptyset$
without loss of generality.
We can further assume that $A$ intersects
the boundary of $l^t(\Delta_0)$ only in the vertical
boundary curves of $l^t(\Delta_0)$, namely
$$
A\cap \partial^hl^t(\Delta_0)=\emptyset.
$$
Indeed, if the equality above is not true
for all $A=A^n$ with $n$ large, then we could prove $\Lambda\cap W^{ss}_X(\sigma)\neq
\{\sigma\}$ by taking limit points
of a sequence $p_n\in A^n\cap \partial^hS^t$
(note that $\partial^h l^t(\Delta_0)\subset
\partial^h S^t$ by definition).

We denote $S=l^t(\Delta_0)$,
$(x,y)=(x^t,y^t)$ and $\Pi=\Pi^t$ for simplicity.
In particular, the line field $F$ is transverse to
$\Pi$.

As $S\cap A$ is compact in $S$,
there is $p\in S\cap A$
satisfying
$$
dist(\Pi(S^t\cap A),0)
=dist(\Pi(p),0),
$$
where $dist$ denotes the distance in $[-\Delta_0,\Delta_0]$.

Now, $p\in A$ and so $W^u_{Z}(p)$
is a well defined two-dimensional
submanifold.
The dominance in Definition \ref{d2}-(2) implies that
$$
T_z(W^u_{Z}(p))=E^{c}_z
$$
for every $z\in W^u_{Z}(p)$, and so,
$$
T_z(W^u_{Z}(p))\cap T_zS=E^c_z\cap T_zS=F_z
$$
for every $z\in W^u_{Z}(p)\cap S$.

As $W^u_{Z}(p)\cap S$ is transversal,
we have that
$W^u_{Z}(p)\cap S$ contains a curve
$C$ whose interior contains $p$ as in Figure
The last equality implies that $C$ is tangent
to $F$.

\begin{figure}[htv] 
\centerline{
\psfig{figure=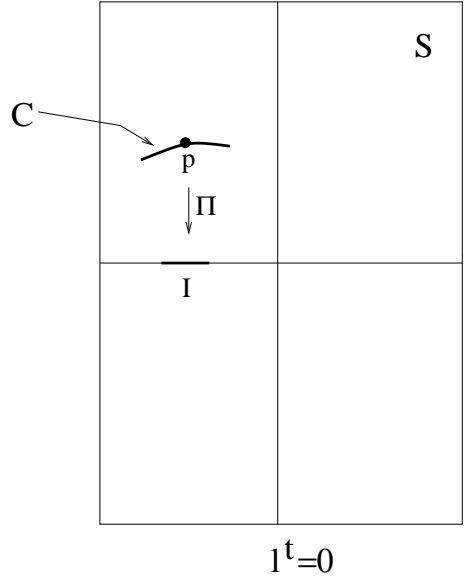,height=3in}}
\caption{\label{f.2} The curve $C$.}
\end{figure}

As $F$ is transverse
to $\Pi$ we have that
$C$ is transverse to $\Pi$ as well.
We conclude that $\Pi(C)$ contains an open
interval $I\subset [-\Delta_0,\Delta_0]$ with
$\Pi(p)\in int(I)$.
So, there
is $z_0\in C$ such that
$$
dist(\Pi(z_0),0)<dist(\Pi(p),0).
$$

Note that $C\subset S\cap A$ since
$A$ is an attractor of $Z$, $p\in A$
and $C\subset W^u_{Z}(p)$.
As $A\cap \partial^hS=\emptyset$
we conclude that
$$
dist(\Pi(S\cap A),0)=0.
$$
As $A$ is closed, this last equality
implies
$$
A\cap l^t\neq\emptyset.
$$
From this it would follow that $\sigma=
\sigma_{n_0}\in A=A^{n_0}$,
a contradiction
since $A$ is a hyperbolic strange attractor of $Z$.
This finishes the proof in Case II.

The proof of Theorem \ref{thB} is completed.

\qed

\section{Proof of the corollaries}

{\flushleft{\bf Proof of Corollary \ref{bound}: }}
Fix the neighborhood $U$ as in Theorem \ref{thB}.
Let
$n$ be the number of singularities
of $X$ in $\Lambda$.
By Theorem \ref{thB}
there is a neighborhood ${\cal W}$ of
$X$ in ${\cal X}^r(M)$
such that every attractor
in $U$ of every $Y\in {\cal W}$ is singular.
Shrinking ${\cal W}$ if necessary
we can assume that
the number of singularities of $Y$ in
$U$ is also $n$.
As the family of attractors of $Y$ in $U$
is pairwise disjoint,
we conclude that every $Y\in {\cal W}$ has at most
$n$ attractors in $U$.
This proves the result.

\qed

{\flushleft{\bf Proof of Corollary \ref{thA}: }}
Let $X$ be a flow in ${\cal X}^r(M)$
exhibiting a
singular hyperbolic attractor $\Lambda$, $r\geq 1$.
Fix the neighborhood $U$ as in Theorem \ref{thB}.
Then every attractor in $U$ of every flow
$C^r$ close to $X$ is singular.

Suppose by contradiction that
the conclusion of the theorem fails.
Then there is a sequence $X^n\to X$ such that
every singularity of $X^n$ is isolated
in the nonwandering set.
Choosing $n$ large
we have that $X^n$ is $C^r$ close to $X$.
In particular, the sequence $X^n$ satisfies
the following property :

\begin{description}
\item{(*).}
Every attractor in $U$ of
$X^n$ is singular.
\end{description}

On the contrary, we claim that every $X^n$ has a hyperbolic attractor
$A^n$ in $U$.
This will follows from the fact that
every singularity of $X^n$ is isolated
in the nonwandering set. Indeed,
define
$$
C^n=
(\Omega(X^n)\cap U)\setminus Sing(X^n).
$$
Then $C^n$ is a compact invariant
set of $X^n$ since each sigularity is isolated
in the nonwandering set. By Lemma \ref{l1}
we have that $C^n$ is hyperbolic
since it has no singularities.
But
$$
\Omega(X^n)\cap U=
C^n\cup \left(Sing(X^n)\cap U\right)
$$
and $Sing(X^n)\cap U$ is a hyperbolic set of $X^n$ too.
Then
$
\Omega(X^n)\cap U
$
is a hyperbolic set of $X^n$.
As $dim(M)=3$ it follows from \cite{np} that
the closed orbits of $X^n$ in $U$ are
dense in $\Omega(X^n)\cap U$ for all $n$
(althought this was proved
in [PN] for surface
diffeomorphisms the same proof works
for three-dimensional flows).
By the Spectral Decomposition Theorem
\cite{pt} we have that $X^n$ has a hyperbolic attractor
$A^n\subset U$ and the claim follows.

By (*) we have that $A^n$ is,
at the same time, a hyperbolic attractor
and a singular hyperbolic attractor
of $X^n$. But this is impossible by
Remark \ref{r3}-(1).
This finishes the proof.

\qed

{\flushleft{\bf Proof of Corollary \ref{c1}: }}
First observe that  no notrivial attractor
can be accumulated by attracting or repelling
closed orbits of the flow.
Then, by applying the arguments
in \cite{mpa}, we can prove that there is
${\cal R}\subset {\cal X}^1(M)$ residual
such that every nontrivial attractor with singularities
of every flow
$X$ in ${\cal R}$ is singular hyperbolic.
Then the conclusion follows applying
Corollary \ref{thA} to $r=1$.

\qed

We finish with an observation.
By using Corollary \ref{Coro}
it is possible to construct
examples
of $C^r$ flows without $C^r$ robust singular
attractors
exhibiting
singularities $C^r$ stably non isolated
in the nonwandering set.
The geometric
model described in \cite[Appendix]{mp}
is one of such examples.

\medskip 

\flushleft
C. A. Morales\\
Instituto de Matem\'atica, Universidade Federal do Rio de Janeiro\\
P. O. Box 68530, 21945--970 Rio de Janeiro, Brazil\\
E-mail: morales@impa.br
\end{document}